\documentclass[10pt,twocolumn]{asme2ej}

\usepackage{graphicx}
\usepackage{hyperref}
\usepackage{cite}

\usepackage{amsmath,amsbsy,amssymb,amsthm}

\newtheorem{theorem}{Theorem}

\newtheorem{remark}{Remark}

\begin{document}

\title{Fractional Order Version of the HJB Equation\thanks{This 
is a preprint of a paper whose final and definite form is with
\emph{Journal of Computational and Nonlinear Dynamics},
ISSN 1555-1415, eISSN 1555-1423, CODEN: JCNDDM. 
Submitted 28-June-2018; Revised 15-Sept-2018; Accepted 28-Oct-2018.}}

\author{Abolhassan Razminia
\affiliation{Ph.D. in Control Systems,\\
Associate Professor,\\
Electrical Engineering Department,\\
School of Engineering,\\
Persian Gulf University,\\
P.~O. Box 75169, Bushehr, Iran\\
Email: razminia@pgu.ac.ir}}

\author{Mehdi AsadiZadehShiraz
\affiliation{M.Sc. in Control Systems,\\
Electronic and Electrical\\ 
Engineering Department,\\
Shiraz University of Technology,\\ 
P.~O. Box 71555-313, Shiraz, Iran\\
Email: M.Asadizadehshiraz@sutech.ac.ir}}

\author{Delfim F. M. Torres\thanks{Corresponding author.}
\affiliation{Ph.D. and D.Sc. (Habilitation) in Mathematics,\\
Full Professor, Coordinator of the R\&D Unit CIDMA,\\
Center for Research and Development in Mathematics and Applications,\\
Department of Mathematics, University of Aveiro,
Aveiro 3810-193, Portugal\\
Email: delfim@ua.pt}}


\date{Received: June 28, 2018 / Revised: Sept 15, 2018 / Accepted: Oct 28, 2018}

\maketitle


\begin{abstract}
{\it We consider an extension of the well-known Hamilton--Jacobi--Bellman (HJB) 
equation for fractional order dynamical systems in which a generalized performance 
index is considered for the related optimal control problem. 
Owing to the nonlocality of the fractional order operators, the classical HJB equation, 
in the usual form, does not hold true for fractional problems. Effectiveness of the proposed 
technique is illustrated through a numerical example.}

\bigskip


\medskip

\end{abstract}


\section{Introduction}

During the last sixty years, dynamic optimization and issues related 
to optimal control theory have received a lot of attention. 
Various theories, and a large number of applications of optimal control, 
can be considered as an indicator of the impact of the theory 
on the science and industry. As a short list of applications, 
we can mention the technology of wave energy converters in an optimal manner \cite{zou}, 
optimal control of gantry cranes \cite{kolar}, emission management in diesel engines \cite{donk}, 
thermic processes \cite{sal}, and epidemiology \cite{MR3721854}.

As is well known, one of the most basic requirements in optimal control theory 
is the modeling of the process. The more accurate the model, the better the control obtained. 
Parallel to advancement in optimal control theory, several mathematical tools have been developed 
to model the process to be controlled. One of these tools is fractional calculus (FC), 
which is an extension of the traditional integer order calculus \cite{kilbas,igor}. 
FC has affected the control engineering discipline in two aspects: getting superior models 
for the processes, and a robust structure of the closed-loop control system \cite{monj,sabb}.
Applications of FC have been explored to various fields of science and engineering, 
including control engineering \cite{salim}, chaotic systems \cite{mee,kumar}, 
reservoir engineering \cite{mee1}, diffusive processes \cite{mee2}, and so on \cite{MR3673710}.
For the design of variable-order fractional proportional--integral--derivative 
controllers for linear dynamical systems, see \cite{MR3787676}.
For applications of fractional calculus on the nutrition of pregnant women 
and the health of newborns and nursing mothers, see \cite{MR3820541};
fractional models of HIV-infection and their potential to extract 
new hidden features of biological complex systems are investigated in \cite{MR3844370}.
In \cite{TS:02}, a time-fractional modified Kawahara equation, describing 
the generation of non-linear water-waves in the long-wavelength regime, is proposed
and studied. A fractional model for convective straight fins with temperature-dependent 
thermal conductivity is found in \cite{TS:01}, while
a fractional model for ion acoustic plasma waves is investigated in \cite{MR3735103}.
In \cite{MR3740466}, a fractional Fitzhugh--Nagumo equation is employed to describe 
the transmission of nerve impulses; in \cite{MR3787704}, 
a fractional model of Lienard's equation is used to describe oscillating circuits.

Generally, there are two main approaches in solving an optimal control problem: 
Bellman's Dynamic Programming and Pontryagin's Maximum Principle (PMP). 
The former presents a necessary and sufficient condition of optimality 
whereas the latter provides necessary conditions. Extending the optimal 
control problem to the fractional order context has been done via PMP 
using different approaches \cite{MR3529374,ag,tor1,MR3673702,me1,tor2,me2}. 
In \cite{MR3798633}, an efficient optimal control scheme is proposed for 
nonlinear fractional-order systems with disturbances,
while optimization of fractional systems with derivatives 
of distributed order is investigated in \cite{Zaky2018}.
In \cite{dz,pr}, the dynamic programming procedure has been extended 
for fractional discrete-time systems. However, to the best of our knowledge, 
the fractional order version of the well-known HJB equation has not been studied thoroughly. 
A controlled continuous time random walk, and their position-dependent extensions, 
have been studied in the framework of fractional calculus \cite{vassili}, 
but their scope and nature are completely different from the study in the current manuscript. 

Here, a general version of the HJB equation is presented. Our main contributions 
are twofold: we define the performance index of the optimal control problem 
in a very general form, with the help of a fractional order operator, and, 
based on the optimality principle, we develop the HJB equation in the fractional context, 
which we denote as the Fr-HJB equation. Since the problem is inherently difficult, 
analytical solutions, even for very simple problems, are in general impossible or
very difficult to obtain. Moreover, numerical simulations for the Fr-HJB equation 
are not an obvious issue. Although different approaches exist for solving the 
classical HJB equation, for the case of the Fr-HJB a reliable numerical 
technique must be chosen, guaranteeing convergence and stability.

The paper is organized as follows. Some preliminaries on FC are given 
in Section~\ref{sec:2}. Section~\ref{sec:3} is devoted to the problem 
statement, in which a formal definition of the optimal control problem 
and some necessary tools are presented. The Hamilton--Jacobi--Bellman 
equation, in the context of fractional order systems, namely the Fr-HJB equation, 
is investigated in Section~\ref{sec:4}. In Section~\ref{sec:5}, 
an optimal control problem is explored via the Fr-HJB equation. Finally, 
Section~\ref{sec:6} concludes the paper with some remarks 
and future research directions.


\section{Preliminaries}
\label{sec:2}

Let $x\in L_1[t_0,t_f]$ be a function of $t\in [t_0,t_f]$. 
The fractional integral of order $v \geq 0$ is defined, 
in the sense of Riemann--Liouville, as follows:
\begin{align}
&\text{Left operator,} 
\quad {}_{t_0}I_t^vx(t)=\frac{1}{\Gamma (v)}
\int_{t_0}^t (t-\tau)^{v-1}x(\tau)\mathrm{d}\tau,\\
&\text{Right operator,} \quad {}_tI_{t_f}^vx(t)
=\frac{1}{\Gamma(v)} \int_t^{t_f}(\tau - t)^{v-1}x(\tau)\mathrm{d}\tau,
\end{align}
where $\Gamma(\cdot)$ is given by
\begin{equation}
\Gamma (v)=\int_0^{\infty} z^{v-1}e^{-z}\mathrm{d}z.
\end{equation}
For $v>0$, we denote the space of functions that can be represented 
by a left (right) RL-integral of order $v$ of some $C([a,b])$-function 
by ${}_aI_t^v([a,b])$ (${}_tI_b^v([a,b])$), where $a<t<b$.

Based on the definition of fractional integrals, the (left) 
fractional derivative operators in the sense of Riemann--Liouville 
and Caputo are defined, respectively, as follows:
\begin{align}
{}_aD_t^q =D^n \circ {}_aI_t^{n-q},\\
{}_a^CD_t^q ={}_aI_t^{n-q} \circ D^n,
\end{align}
where $n \in \mathbb{Z}^+$ and $n-1 \leq q <n$. 
The relationship between Caputo and Riemann--Liouville derivatives 
is given by the following formula:
\begin{equation}
\label{rlc}
{}_a^CD_t^qx(t)={}_aD_t^qx(t)-\sum_{k=0}^{n-1} 
\frac{x^{(k)}(a)}{\Gamma (k-q+1)}(t-a)^{k-q},
\end{equation}
where $n-1 \leq q < n$. Along the work, 
we consider $0<q<1$. In this case, the above relation reduces to
\begin{equation}
{}_a^CD_t^qx(t)={}_aD_t^qx(t)-\frac{x(a)}{\Gamma (1-q)}(t-a)^{-q},
\end{equation}
which implies that the difference between Caputo and Riemann--Liouville 
derivatives depends on the initial value of $x(t)$. In addition, 
as stated in \cite{li}, for $0<q<1$ we have, in the limit,
\begin{equation}\label{dcd}
\lim_{q\rightarrow 1^-}{}_a^CD_t ^qx(t)
=\lim_{q\rightarrow 1^-}{}_aD_t ^qx(t)=\dot{x}(t).
\end{equation}
Therefore, the Caputo and Riemann--Liouville 
derivatives are consistent with
the standard integer order derivatives.

One of the important relationships between fractional order 
derivatives and the integer ones can be stated as in \cite{Atana}:
\begin{multline}
\label{tag}
{}_aD_t^qx(t)=A(q)(t-a)^{-q}x(t)+B(q)(t-a)^{1-q}\dot{x}(t)\\
-\sum_{p=2}^{\infty}C(q,p)(t-a)^{1-p-q}W_p(t),
\end{multline}
where for $p=2,3,\ldots$ function $W_p(t)$ solves 
the initial value problem
\begin{equation}
\label{eq:Wp}
\begin{cases}
\dot{W}_p(t)=(1-p)(t-a)^{p-2}x(t),\\
W_p(a)=0,
\end{cases}
\end{equation}
and the coefficients $A(q),B(q)$ and $C(q,p)$ 
are determined by the following formulas:
\begin{align}
A(q)
&= \frac{1}{\Gamma (1-q)}\left( 1+\sum_{p=2}^{\infty} 
\frac{\Gamma (p-1+q)}{\Gamma (q)(p-1)!} \right),\label{eq:A}\\
B(q)&= \frac{1}{\Gamma (2-q)}\left( 1+\sum_{p=1}^{\infty} 
\frac{\Gamma (p-1+q)}{\Gamma (q)(p-1)!} \right),\label{eq:B}\\
C(q,p)&= \frac{1}{\Gamma (2-q)\Gamma (q-1)} 
\cdot \frac{\Gamma (p-1+q)}{(p-1)!}. \label{eq:C}
\end{align}
The backward compatibility of relation \eqref{tag} for the case $q\rightarrow 1$ 
can be easily proven by considering the properties of the gamma function.


\section{Problem statement}
\label{sec:3}

Consider a plant with the dynamical control system
\begin{equation}\label{plant}
{}_{t_0}^CD_t^{\mathbf{q}} \mathbf{x}(t)=\mathbf{f}(t,\mathbf{x}(t),\mathbf{u}(t)),
\end{equation}
where $\mathbf{q}=\left[
\begin{array}{cccc}
q_1 & q_2 & \ldots & q_n
\end{array}
\right]
$ is the order of differentiation in the sense of Caputo, so that 
$0< q_i <1$ for $i=1,2,\ldots , n$, $\mathbf{f}$ is a smooth vector field 
for the pseudo-states $\mathbf{x}(t)\in \mathcal{X} \subseteq \mathbb{R}^n$ 
and $\mathbf{u}(t) \in \mathcal{U} \subset \mathbb{R}^m$ 
is the control input vector, where $\mathcal{U}$ is a compact set. 
Moreover, the initial state is assumed to be known as 
$\mathbf{x}(t_0)=\mathbf{x}_0$, where the initial time $t=t_0$ is the starting point 
of the dynamic system, which implicitly assumes that all information from 
$-\infty$ to $t_0$ is summarized in $\mathbf{x}_0$.

Let $\mathfrak{O}=
\left[
\begin{array}{cccc}
\mathcal{O}_1 & \mathcal{O}_2 & \ldots & \mathcal{O}_r
\end{array}
\right]$ be a vector-operator and $\mathbf{g}=\left[
\begin{array}{cccc}
g_1 & g_2 & \ldots & g_r
\end{array}
\right]
$ be a vector function so that the operator $\mathcal{O}_j$ affects the operand $g_j$. 
Then, the generalized dot-product is defined as follows:
\begin{equation}
\mathfrak{O} \odot \mathbf{g}=\sum_{j=1}^r  \mathcal{O}_j g_j,
\end{equation}
where $\mathcal{O}_jg_j \in \mathbb{R}$.
Let us define the vector operator as
\begin{equation}
{}_{t_0}\mathbf{I}_{t_f}^\mathbf{v}:=\left[
\begin{array}{cccc}
{}_{t_0}I_{t_f}^{v_1} & {}_{t_0}I_{t_f}^{v_2} & \ldots & {}_{t_0}I_{t_f}^{v_r}
\end{array}
\right],
\end{equation}
where $0\leq v_i\leq 2$ for $i=1,2,\ldots , r$ is called the tuning factor 
of the performance index. Here, we study the fixed-final-time problem 
in which $t_f$ is specified beforehand. Therefore, the Riemann--Liouville integral 
of left type is preferable than the right one:
\begin{equation}
{}_{t_0}I_{t_f}^v g(t)=\frac{1}{\Gamma (v)}
\int_{t_0}^{t_f}(t_f- \tau)^{v-1}g(\tau)\mathrm{d}\tau.
\end{equation}
The main motivation for introducing such performance index can be explained 
by noticing the kernel. Depending on the nature of the problem, 
we can choose the order of integration, $v$, so that the desired criterion is achieved. 
Indeed, there are two highlighted cases:
\begin{description}
\item
\textbf{Expensive initial behavior}, when the initial behavior is more 
important than the final ones, so that we can select the order of the 
integral greater than one: $1< v_j<2$. By this choice, the kernel 
$(t_f- \tau)^{v-1}$ is a larger quantity for initial times, 
whereas for the final times it will be smaller.

\item \textbf{Cheap initial behavior}, for which, in contrast, the final behavior 
is more important, so that we can select the tuning factor of the performance index 
$v_j$ as $0<v_j<1$. For this case, a bigger weight is imposed on the final time 
behavior and smaller weights on the initial times.
\end{description}
Apparently, the tuning factor $\mathbf{v}$ weights the integrand naturally in the time basis. 
Such time-filtering cannot be so easily included in the classical performance index. 
Moreover, since a vector type of order $\mathbf{v}$ has been chosen, both 
expensive and cheap cases can be considered compactly in an unified expression. 
Moreover, under some mild assumptions, we can generalize the performance index 
for $v_j<0$, which implies that, using a unified framework, we can consider 
a derivative cost functional:
\begin{equation}
{}_{t_0}I_t^{v_j} \rightarrow {}_{t_0}D_t^{-v_j}.
\end{equation}
Such issues may appear in limited-saturation control problems 
in which the derivative of the control signal has to be constrained. 
In addition to these cases, a combined expression can be constituted, 
in which some integrals, with different orders, depending on the nature of the problem,
and some derivatives with different orders, appear. Therefore, the proposed generalized 
performance index provides a general format for evaluating the optimality of a dynamical system.

Based on this short motivation to the use of a generalized performance index, 
we formulate an optimum behavior of the given plant by minimizing 
the following cost functional:
\begin{equation}
\label{pi}
J={}_{t_0}\mathbf{I}_{t_f}^\mathbf{v}\odot \mathbf{g},
\end{equation}
where the minimization is taken over the control signal $\mathbf{u}$,
assumed to be measurable.
The operand $g_j$ can be considered as the modified Lagrangian 
for the dynamical system, which is in general a mapping 
$g_j(t,\mathbf{x},\mathbf{u}): \mathbb{R}^+\times \mathbb{R}^n 
\times \mathbb{R}^m\mapsto \mathbb{R}$ when its related operator 
$\mathcal{O}_j = {}_{t_0}I_{t_f}^{v_j}$ 
has nonzero order, $v_j\neq 0$. For the operator 
of order zero (nonintegral part), the operand is assumed to be 
$g_j(t_f,\mathbf{x}(t_f))$.

\begin{remark}
It can be easily seen that the classical Bolza problem of optimal control theory 
is a special class of the above generalized problem in which 
$r=2$, $v_1=0$, $v_2=1$, $g_1=h(t_f,\mathbf{x}(t_f))$, 
and $g_2=g(t,\mathbf{x}(t),\mathbf{u}(t))$. In this case, 
the performance index is in the Bolza form \cite{kirk}:
\begin{equation}
J=h(t_f,\mathbf{x}(t_f))+\int_{t_0}^{t_f} g(t,\mathbf{x}(t),\mathbf{u}(t)) \mathrm{d}t.
\end{equation}
As can be seen, the fractional order performance index \eqref{pi}
has the backward compatibility property, i.e., considering integer order parameters, 
the classical optimal control problem is obtained.
\end{remark}

An action-like integral for problems of the calculus of variations  
has been introduced in \cite{nabul1,nabul2,nabul3}. However, 
our main motivation here to define such generalized performance index 
is completely different and relies, basically, on the natural 
weighting and backward compatibility properties.


\section{HJB equation: fractional order version}
\label{sec:4}

Consider a generalized performance index over the interval $[t,t_f]$, 
$t\leq t_f$, of any control sequence $\mathbf{u}(\tau), t\leq \tau \leq t_f$:
\begin{equation}
J(t,\mathbf{x}(t),\mathbf{u}(\tau))={}_{t}\mathbf{I}_{t_f}^\mathbf{v}\odot \mathbf{g}.
\end{equation}
Based on Bellman's optimality principle \cite{kirk}, the first step 
is computing the optimal cost at the final time $t_f$. This value can 
be obtained by observing the non-integral terms of $J$, i.e., 
the terms whose integration order is zero: $v_k=0$. 
Let us denote these indexes by the set $\mathcal{K}$:
\begin{equation}
\mathcal{K}=\{ k: v_k=0, k=1,2,\ldots , r \}.
\end{equation}
For such terms, the optimum values of $V$ at the final time are
\begin{equation}
\label{fin}
V(t_f,\mathbf{x}(t_f))=\sum_{k\in \mathcal{K}}g_k.
\end{equation}
For the empty set $\mathcal{K}$, we consider $V(t_f,\mathbf{x}(t_f))=0$ 
as the boundary value of the $V$-function. Now we use the backward trend 
of the dynamic programming procedure. Clearly, the goal is to pick 
$\mathbf{u}(\tau), t\leq \tau \leq t_f$, to minimize the cost functional
\begin{equation}
J^*(t,\mathbf{x}(t))=V(t,\mathbf{x})
=\inf_{\substack{\mathbf{u}(\tau) \in U,\\ 
t\leq \tau \leq t_f}} J(t,\mathbf{x}(t),\mathbf{u}(\tau)).
\end{equation}
By writing the value function explicitly, and then splitting 
the control horizon $[t,t_f]$ into two subinterval $[t,t+\Delta t]$ 
and $[t+\Delta t, t_f]$, we get:
\begin{align}
\label{v1}
V(t,\mathbf{x}(t))
&=\inf _{\substack{\mathbf{u}(\tau)\in U\\ 
t\leq \tau \leq t_f}} \left\lbrace  {}_t\mathbf{I}_{t_f}^{\mathbf{v}}
\odot \mathbf{g}   \right\rbrace \nonumber\\
&= \inf_{\substack{\mathbf{u}(\tau)\in U\\ 
t\leq \tau \leq t_f}} \left\lbrace  {}_t\mathbf{I}_{(t+\Delta t)}^{\mathbf{v}} 
\odot \mathbf{g} + {}_{(t+\Delta t)}\mathbf{I}_{t_f}^{\mathbf{v}} 
\odot \mathbf{g} \right\rbrace.
\end{align}
In this stage, at time $t+\Delta t$, the system will be implicitly 
at pseudo-state $\mathbf{x}(t+\Delta t)$. But from the principle of optimality, 
we can write the optimal cost-to-go from this state as
\[
V(t+\Delta t,\mathbf{x}(t+\Delta t)).
\]
Thus, we can rewrite the cost calculation in Eq. \eqref{v1} as:
\begin{align}
\label{baz}
V(t,\mathbf{x}(t))
&= \inf _{\substack{\mathbf{u}(\tau)\in U\\ 
t\leq \tau \leq t_f}} \left\lbrace  {}_t\mathbf{I}_{(t+\Delta t)}^{\mathbf{v}} 
\odot \mathbf{g} +V(t+\Delta t,\mathbf{x}(t+\Delta t)) \right\rbrace \nonumber\\
&\approx \inf _{\substack{\mathbf{u}(\tau)\in U\\ 
t\leq \tau \leq t_f}} \Bigg\lbrace \sum_{j=1}^r 
\frac{1}{\Gamma (v_j)} (t_f-t)^{v_j-1}g_j(t,
\mathbf{x}(t),\mathbf{u}(t))\Delta t\\
& \qquad \qquad \quad 
+V(t+\Delta t,\mathbf{x}(t+\Delta t))\Bigg\rbrace.\nonumber
\end{align}
Assuming that $V$ has bounded second derivatives 
in both arguments, one can expand this cost as a Taylor series about $(t,\mathbf{x}(t))$:
\begin{align}\label{vv}
V&(t+\Delta t,\mathbf{x}(t+\Delta t))
\simeq V(t,\mathbf{x}(t))
+\left[\frac{\partial V}{\partial t}(t,\mathbf{x}(t)) \right]\Delta t \nonumber\\
&+\left[\frac{\partial V}{\partial \mathbf{x}}(t,\mathbf{x}(t)) \right]^T(\mathbf{x}(t+\Delta t)
-\mathbf{x}(t))\nonumber \\
&\approx V(t,\mathbf{x}(t))+V_t(t,\mathbf{x}(t))\Delta t 
+ V_{\mathbf{x}}^T(t,\mathbf{x}(t)) \dot{\mathbf{x}}\Delta t,
\end{align}
in which for small $\Delta t$, the term $(\mathbf{x}(t+\Delta t)-\mathbf{x}(t))$ 
has been replaced by $\dot{\mathbf{x}}\Delta t$. Now, resorting to Eqs. \eqref{tag} 
and \eqref{rlc}, we can replace $\dot{\mathbf{x}}$ 
by the following relation (component wise):
\begin{align*}
\dot{x}_i
&=\frac{{}_aD_t^{q_i}x_i(t)-A(q_i)(t-a)^{-q_i}x_i(t)}{B(q_i)(t-a)^{1-q_i}}\\
&\qquad +\frac{\sum_{p=2}^{\infty} C(q_i,p)(t-a)^{1-p-q_i}W_p(t)}{B(q_i)(t-a)^{1-q_i}}\\
&=\frac{{}_a^CD_t^{q_i}x_i(t) -k_i(t,q_i,x_i)}{B(q_i)(t-a)^{1-q_i}},
\end{align*}
where
\begin{multline}
k_i(t,q_i,x_i)=-\frac{x_i(a)}{\Gamma (1-q_i)}(t-a)^{-q_i}
+ A(q_i)(t-a)^{-q_i}x_i(t)\\
-\sum_{p=2}^{\infty} C(q_i,p)(t-a)^{1-p-q_i}W_p(t)
\end{multline}
and $A(\cdot)$, $B(\cdot)$, $C(\cdot)$ and $W_p(\cdot)$ 
are defined in \eqref{eq:Wp}--\eqref{eq:C}. Let
\begin{equation}
\label{f}
\tilde{\mathbf{f}}=\left[
\begin{array}{cccc}
\frac{f_1-k_1(t,q_1,x_1)}{B(q_1)(t-a)^{1-q_1}}&
\frac{f_2-k_2(t,q_2,x_2)}{B(q_2)(t-a)^{1-q_2}}&
\ldots&
\frac{f_n-k_n(t,q_n,x_n)}{B(q_n)(t-a)^{1-q_n}}
\end{array}
\right]^T.
\end{equation}
By this transformation, the equivalent equation describing the system is 
$\dot{\mathbf{x}}=\tilde{\mathbf{f}}(t,\mathbf{x},\mathbf{u})$.
Therefore, Eq.~\eqref{vv} reduces to
\begin{equation*}
V(t+\Delta t, \mathbf{x}(t+\Delta t)) \approx V(t,\mathbf{x}(t))
+V_t(t,\mathbf{x}(t))\Delta t + V_{\mathbf{x}}^T(t,\mathbf{x}(t)) 
\tilde{\mathbf{f}}\Delta t.
\end{equation*}
Thus, we can simplify Eq. \eqref{baz} in the following form:
\begin{align*}
V(t,\mathbf{x}(t))\approx
\inf _{\substack{\mathbf{u}(\tau)\in \mathcal{U}\\ 
t\leq \tau \leq t_f}} \Bigg\lbrace \sum_{j=1}^r 
\frac{1}{\Gamma (v_j)} (t_f-t)^{v_j-1}g_j(t,\mathbf{x}(t),\mathbf{u}(t))\Delta t\\
+V(t,\mathbf{x}(t))+V_t(t,\mathbf{x}(t))\Delta t 
+ V_{\mathbf{x}}^T(t,\mathbf{x}(t)) \tilde{\mathbf{f}}\Delta t \Bigg\rbrace.
\end{align*}
Since the minimization is taken over $\mathbf{u}(\cdot)$, 
the term $V(t,\mathbf{x}(t))$ can be canceled from both sides. 
The result just proved is summarized in the following theorem.

\begin{theorem}[Fractional HJB equation]
\label{thm:FHJB}
Consider the plant described by Eq. \eqref{plant} and the performance index 
\eqref{pi}. Assume that $(\mathbf{u},\mathbf{x})$ is the optimal pair, 
which minimizes the performance index $J$. Then the value-function 
$V(t,\mathbf{x})$ satisfies
\begin{multline}
-V_t(t,\mathbf{x})=\min _{\substack{\mathbf{u}(\tau)\in U\\ 
t\leq \tau \leq t_f}} \Bigg\lbrace \sum_{j=1}^r 
\frac{1}{\Gamma (v_j)} (t_f-t)^{v_j-1}g_j(t,\mathbf{x}(t),\mathbf{u}(t))\\
+ V_{\mathbf{x}}^T(t,\mathbf{x}(t)) \tilde{\mathbf{f}} \Bigg\rbrace,
\end{multline}
where
$\tilde{\mathbf{f}}$ is defined in Eq. \eqref{f} and the boundary value 
of $V$ is set as Eq. \eqref{fin}. Moreover, the optimal cost of the system 
is given by
\begin{equation}
J^*=V(t_0,\mathbf{x}_0).
\end{equation}
\end{theorem}

\begin{remark}
Theorem~\ref{thm:FHJB} is backward compatible.
Indeed, it is easy to see that for the limit case where 
$q_i=1$, $r=2$, $v_1=0$, $v_2=1$, $g_1=h(t_f,\mathbf{x}(t_f))$
and $g_2=g(t,\mathbf{x},\mathbf{u})$, the result reduces 
to the classical HJB equation. In this reduction, we can resort  
to Eq. \eqref{dcd}.
\end{remark}


\section{Discussion}
\label{sec:5}

In this section some discussions about the considered problem are presented. 
Several examples for the simplest case $v=r=1$ have been investigated in \cite{seyed}, 
wherein the dynamic of the system has been assumed to be described by 
the Riemann--Liouville derivative. Here, we formulate a more general problem 
and consider the Fr-HJB equation developed in the previous two sections.

Consider a plant with the following dynamics described by the Caputo derivative:
\begin{align}
{}_0^CD_t^{0.2}x_1(t)&=x_2(t)+u(t),\label{ex:cs1}\\
{}_0^CD_t^{0.7}x_2(t)&=-x_1(t),\label{ex:cs2}
\end{align}
subject to the initial conditions $x_1(0)=1$ and $x_2(0)=0.5$.
It is desired to find a control $u$ such that the following 
performance index is minimized:
\begin{equation}
\label{ex:J}
J={}_0I_1^{0.3} \left( x_1^2(t)+x_2^2(t) \Big)+ {}_0I_1^{0.4}\Big(x_1^2(t)+ u^2(t) \right).
\end{equation}
As can be seen, the performance index is free of nonintegral terms. 
Thus, the final value of the $V$-function is set as zero:
\begin{equation}
V(1,x_1(1),x_2(1))=0.
\end{equation}
In this case, the fractional HJB equation is given by
\begin{align*}
-\frac{\partial V}{\partial t}(t,\mathbf{x})
=\min_{\substack{\mathbf{u}(\tau)\in U\\ 
t\leq \tau \leq 1}} \Bigg\{ \frac{1}{\Gamma (0.3)}
 (1-t)^{-0.7}(x_1^2(t)+x_2^2(t))\\
+\frac{1}{\Gamma (0.4)} (1-t)^{-0.6}(x_1^2(t)+u^2(t))
+ \left(\frac{\partial V}{\partial 
\mathbf{x}}(t,\mathbf{x}(t))\right)^T\tilde{\mathbf{f}}\Bigg\}.
\end{align*}
Let us constitute the components of the fractional HJB equation:
\begin{align*}
\tilde{\mathbf{f}}
&=\left[
\begin{array}{cc}
\frac{x_2(t)+u(t)-k_1(t,0.2, x_1)}{B(0.2)t^{0.8}} 
& \frac{-x_1(t)-k_2(t,0.7,x_2)}{B(0.7)t^{0.3}}
\end{array}
\right]^T,\\
k_1(t,0.2,x_1)
&=-\frac{t^{-0.2}}{\Gamma (0.8)}+A(0.2)t^{-0.2} x_1(t)\\
&\quad -\sum_{p=2}^{\infty} C(0.2,p)t^{0.8-p}W_{p1}(t),\\
k_2(t,0.7,x_2) &=-\frac{t^{-0.7}}{2\Gamma (0.3)}+ A(0.7)t^{-0.7} x_2(t)\\
&\quad -\sum_{p=2}^{\infty} C(0.7,p)t^{0.3-p}W_{p2}(t),
\end{align*}
where
\begin{equation*}
\begin{cases}
\dot{W}_{pi}(t)=(1-p)t^{p-2}x_i(t),\\
W_{pi}(0)=0.
\end{cases}
\end{equation*}
In the numerical implementation, the infinite series in $A(\cdot)$, 
$B(\cdot)$ and $C(\cdot,\cdot)$ are truncated up to $N$ steps \cite{MR3443073}:
\begin{align}
A(q)&= \frac{1}{\Gamma (1-q)}\left( 
1+\sum_{p=2}^{N_A} \frac{\Gamma (p-1+q)}{\Gamma (q)(p-1)!}\right),\\
B(q)&= \frac{1}{\Gamma (2-q)}\left( 1+\sum_{p=1}^{N_B} 
\frac{\Gamma (p-1+q)}{\Gamma (q)(p-1)!} \right),\\
C(q,p)&= \frac{1}{\Gamma (2-q)\Gamma (q-1)} \cdot \frac{\Gamma (p-1+q)}{(p-1)!}.
\end{align}
In such case, we can explicitly denote the first two coefficients by $A(\cdot,N_A)$ 
and $B(\cdot,N_B)$, where $N_A$ and $N_B$ are the upper bounds of the summations. 
For some significant works for numerical solutions, 
see \cite{MR3443073,yong,sun,xu,chenn,zhoo}.

\begin{figure}
\centering
\includegraphics[scale=0.6]{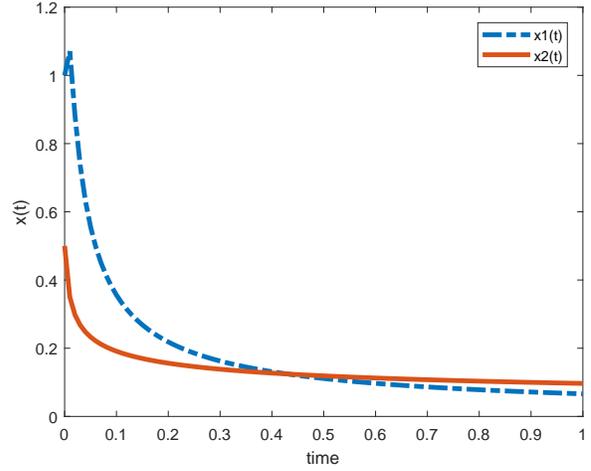}
\caption{Optimal system states of problem \eqref{ex:cs1}--\eqref{ex:J}, 
which converge to the equilibrium points.}\label{Pic1}
\end{figure}

\begin{figure}
\centering
\includegraphics[scale=0.6]{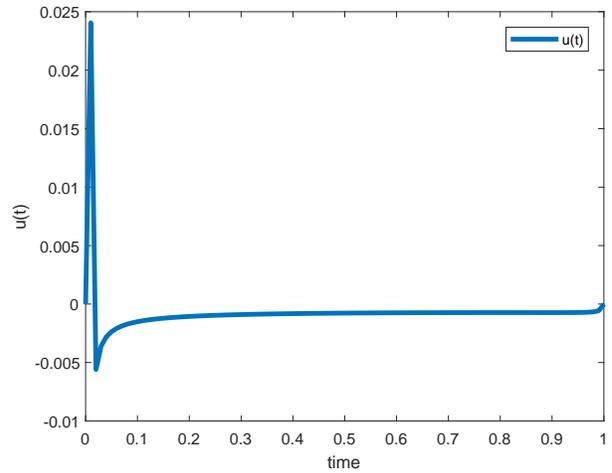}
\caption{The optimal controller of the system \eqref{ex:cs1}--\eqref{ex:cs2}
that minimizes the performance index \eqref{ex:J}.}\label{Pic2}
\end{figure}

It can be clearly seen that some equations must be solved forward while others backward. 
Therefore, in general, these type of problems cannot be solved using conventional 
numerical methods. One of the reliable numerical techniques in solving such optimal 
control problems is the Forward-Backward Sweep Method (FBSM), based on the following five steps:
\begin{enumerate}
\item Guess the initial conditions for the controller $u(t)$ and save it. 
In our example, we have considered $\Delta t=0.01$ and $u(t) = 5$.

\item Acquire and save states $x(t)$ in forward time, based on the given initial conditions 
of the states and the $u(t)$ stored in step one. We have considered $N_A=N_B=10^9$ and 
limited $p$ to $150$ for calculating $W_{pi}$, which results 
in $x_1(1)\simeq0.138$ and $x_2(1)\simeq0.097$. 

\item Obtain and save the co-states in backward, according to their final conditions.
By using the backward path, the vector $V(t)$ can be obtained. For our considered 
initial time we obtain $V(0)\simeq8.3$.

\item Update $u(t)$ with respect to states and co-states 
obtained from steps 2 and 3.

\item Check the variables values and their error rates. 
If the error is small enough, the obtained values are considered valid 
and the process is finished. Otherwise, jump and start from step 2.
\end{enumerate}
Note that the presented method only applies to numerical solutions of problems 
in which the initial condition of $x(t)$ is constant and at the other times are free. 
For the considered problem, the error is defined as follows:
\begin{multline}
error(t)= \sum_{j=1}^r \frac{1}{\Gamma (v_j)} (t_f-t)^{v_j-1}
g_j(t,\mathbf{x}(t),\mathbf{u}(t))\\
+ V_{\mathbf{x}}^T(t,\mathbf{x}(t)) \tilde{\mathbf{f}}+V_t(t,\mathbf{x})
\end{multline}
and
\begin{equation}
 Error=\left(\sum_{n=0}^{\frac{t_f}{\Delta t}} error^2(n \Delta t)\right)^{(1/2)}.
\end{equation}
By using the proposed method, the optimal cost of the system is 
$V(0,x(0))\simeq0.0053$, $x_1(1)\simeq0.0667$, $x_2(1)\simeq0.0970$ 
and $Error\simeq1.06\times10^{-15}$. The state trajectories and the 
controller output signals are shown in Figures~\ref{Pic1} and \ref{Pic2}, 
respectively.


\section{Concluding remarks}
\label{sec:6}

Based on fractional calculus theory, a generalized performance index 
has been defined for a typical optimal control problem in which 
the dynamical system is also of fractional order. We observed that 
the mentioned performance index has two important properties: 
backward compatibility to the integer order case and a 
natural weighting process. Besides these main features, 
one can include, under some mild assumptions, a derivative-type 
performance index in the proposed unified framework. Subsequently, 
based on the optimality principle, we investigated a general optimal 
control problem with a vectorized order performance index.
Thanks to the continuous time dynamic programming theory
and the series expansions for fractional calculus proposed by
\cite{Atana} and further explored in \cite{MR3443073}, a fractional order 
version of the well-known Hamilton--Jacobi--Bellman (HJB) 
equation has been derived (Theorem~\ref{thm:FHJB}). 
Finally, some discussions about the computational 
difficulties of the problems were presented, 
which can be considered as an important future line of research.


\begin{acknowledgment}
Razminia would like to thank Igor Podlubny and Ivo Petras 
for their encouragement and helpful discussions during his visit
at Technical University of Kosice, Kosice, Slovak Republic.
Torres was supported by FCT through CIDMA,
project UID/MAT/04106/2013, and TOCCATA,
project PTDC/EEI-AUT/2933/2014,
funded by FEDER and COMPETE 2020.
The authors are very grateful to two anonymous referees 
for reading their paper carefully and for all 
the constructive remarks and suggestions.
\end{acknowledgment}



\end{document}